\documentclass[11pt]{article}
\usepackage{bbm}
\usepackage{amsfonts}
\usepackage{amsmath}
\usepackage{latexsym,amssymb,amsmath,amsfonts,amsthm}
\makeatletter
\newif\ifmsbmloaded@
\def\loadmsbm{\msbmloaded@true
  \font\tenmsb=msbm10 scaled 1\@ptsize00
  \font\sevenmsb=msbm7 scaled 1\@ptsize00
  \font\fivemsb=msbm5 scaled 1\@ptsize00
  \alloc@8\fam\chardef\sixt@@n\msbfam
  \textfont\msbfam=\tenmsb
  \scriptfont\msbfam=\sevenmsb
  \scriptscriptfont\msbfam=\fivemsb
  }
\def\nonmatherr@#1{\errmessage%
{LateX error: \string#1\space allowed only in math mode}}
\def\Bbb{\relax\ifmmode\expandafter\Bbb@\else
  \expandafter\nonmatherr@\expandafter\Bbb\fi}
\def\Bbb@#1{{\Bbb@@{#1}}}
\def\Bbb@@#1{\fam\msbfam\relax#1}
\@addtoreset{equation}{section}

\makeatother \loadmsbm
\def\R{\Bbb R}

\def\N{\Bbb N}
\def\Z{\Bbb Z}

\def\no{\noindent}

\def\f#1#2{\frac{#1}{#2}}
\def\f{\frac}

\def\pa{\partial}
\def\p{\partial}

\def\na{\nabla}
\def\la{\lambda}

\def\sqr#1#2{\vbox{\hrule height .#2pt
\hbox{\vrule width .#2pt height #1pt \kern #1pt \vrule width
.#2pt}\hrule height .#2pt }}

\def\endproof{\hphantom{MM}\hfill\llap{$\blacksquare$}\goodbreak}
\def\dv{\mbox{\rm div}}

\def\eqdefa{\buildrel\hbox{\tiny def}\over {=\!\!\!=}}

\newcommand{\cB}{{\mathcal B}}
\newcommand{\cC}{{\mathcal C}}
\newcommand{\cD}{{\mathcal D}}

\newcommand{\cF}{{\mathcal F}}
\newcommand{\cG}{{\mathcal G}}

\newcommand{\cS}{{\mathcal S}}

\newcommand{\calZ}{{\mathcal Z}}

\newcommand{\beq}{\begin{equation}}
\newcommand{\eeq}{\end{equation}}
\newcommand{\ben}{\begin{eqnarray}}
\newcommand{\een}{\end{eqnarray}}
\newcommand{\beno}{\begin{eqnarray*}}
\newcommand{\eeno}{\end{eqnarray*}}
\newtheorem{Theorem}{Theorem}[section]
\newtheorem{Lemma}[Theorem]{Lemma}
\newtheorem{Def}{Definition}[section]
\newtheorem{Remark}{Remark}[section]
\newtheorem{Proposition}{Proposition}[section]

\begin{document}
\title{Global well-posedness for the 3D rotating Navier-Stokes equations with highly oscillating initial data}
\author{Qionglei Chen$ ^\dag$  Changxing Miao$ ^\dag$ and Zhifei Zhang$ ^\ddag$
\\[2mm]
{\small $ ^\dag$ Institute of Applied Physics and Computational Mathematics,
Beijing 100088, China}\\
{\small E-mail: chen\_qionglei@iapcm.ac.cn  and
miao\_changxing@iapcm.ac.cn}\\
{\small $ ^\ddag$ School of Mathematical Sciences, Peking University, 100871, P. R. China}\\
{\small E-mail: zfzhang@math.pku.edu.cn}}
\date{2010, March 30}
\maketitle
\vspace{-1.2in} \vspace{.9in} \vspace{0.2cm}

\begin{abstract}
In this paper, we prove the global well-posedness for the 3D rotating Navier-Stokes equations
in the critical functional framework. Especially, this result allows to construct global solutions
for a class of highly oscillating initial data.
\end{abstract}

\section{Introduction}
In this paper, we study the 3D rotating Navier-Stokes equations
\begin{equation}\label{eq:NSrotating}
\left\{
\begin{array}{ll}
u_t-\nu\Delta u+\Omega e_3\times u+u\cdot\nabla u+\nabla p=0,\\
\dv u=0, \\
u(0,x)=u_0(x) ,
\end{array}
\right.
\end{equation}
where $\nu$ denotes the viscosity coefficient of the fluid, $\Omega$
the speed of rotation, $e_3$ the unit vector in $x_3$ direction and
$\Omega e_3\times u$ the Coriolis force. We refer to \cite{CDGG2,Majda,Ped} for  its background in geophysical fluid
dynamics. If the  Coriolis force is neglected,  the equations (\ref{eq:NSrotating}) become the
classical 3D incompressible Navier-Stokes equations
\begin{equation}\label{eq:NS}
\left\{
\begin{array}{ll}
u_t-\nu\Delta u+u\cdot\nabla u+\nabla p=0,\\
\dv u=0, \\
u(0,x)=u_0(x).
\end{array}
\right.
\end{equation}

The global existence of weak solution of (\ref{eq:NSrotating}) can be proved by the classical
compactness method, since we still have the energy estimate
\beno
\|u(t)\|_{L^2}^2+2\int_0^t\|\na u(s)\|^2_{L^2}ds\le \|u_0\|_{L^2}^2.
\eeno
As in the 3D Navier-Stokes equations,
the uniqueness and regularity of weak solutions are also open.
Recently, Giga et al.\cite{Giga1,Giga2,Giga3} studied the local existence of mild solution for a class of
nondecaying initial data which includes a class of almost periodic functions,
and global existence for small data. On the other hand, when the speed $\Omega$ of rotation is fast enough,
the global existence of smooth solution was proved by \cite{BMN1,BMN2,CDGG1,CDGG2}.

For the 3D Navier-Stokes equations,
Fujita and Kato\cite{Fuj-Kat} and Kato\cite{Kato} proved  the local wellposedness for large initial data and
the global wellposedness for small initial data in the homogeneous Sobolev space
$\dot H^{\f 12}$ and the Lebesgue space $L^3$ repectively. These spaces
are all the critical ones, which are revelent to the scaling of the Navier-Stokes equations: if
$(u,p)$ solves  (\ref{eq:NS}), then
\begin{align}\label{scaling}
(u_\lambda(t,x),p_\lambda(t,x))\eqdefa (\lambda u(\lambda^2 t,\lambda x),\lambda^2 p(\lambda^2 t,\lambda x))
\end{align}
is also a solution of (\ref{eq:NS}). The so-called {\bf critical space} is the one such that
the associated norm is invariant under the scaling of \eqref{scaling}.
Recently, Cannone \cite{Can1}(see also \cite{Can-book, Can2,CMP})
generalized it to Besov spaces with negative index of regularity.  More precisely,
he showed that if the initial data satisfies
$$
\|u_0\|_{\dot B^{-1+\f 3p}_{p,\infty}}\le c,\quad p>3
$$
for some small constant $c$, then the Navier-Stokes equations (\ref{eq:NS}) is globally well-posed.
Let us emphasize that this result allows to construct global solutions for highly oscillating initial data
which may have a large norm in $\dot H^{\f 12}$ or $L^3$. A typical example is
\begin{align*}
u_0(x)=\sin\bigl(\f {x_3} {\varepsilon}\bigr)(-\p_2\phi(x), \p_1\phi(x),0)
\end{align*}
where $\phi\in \cS(\R^3)$ and $\varepsilon>0$ is small enough. We
refer to \cite{CG,CZ, CMZ} for some relevant results.  A natural
question is then to prove a theorem of this type for the rotating
Navier-Stokes equations.

 As we know that  Kato's method heavily relies on
the uniform boundedness of Stokes semigroup in $L^p$ and global
$L^p-L^q$ estimates, but Stokes-Coriolis semigroup is not a
uniformly bounded in $L^p$ for $p\not=2$, see Theorem 5 and Theorem
6 in \cite{OSA}.  Standard techniques allow us to prove these
estimates only locally for the Stokes-Coriolis semigroup,  hence one
can obtain the local existence of mild solution in $L^3$ by Kato's
method. Whether one can extend this solution to a global one for
small data in $L^3$ is a very interesting problem.

Very recently, based on the global $L^p-L^q$ estimates with $q\le2\le p$
and $L^q-H^{\f12}$ estimates with $q>3$ for the Stokes-Coriolis semigroup, Hieber and Shibata\cite{HS}
proved the following global result for small data in $H^\f12$.

\begin{Theorem}\label{Th:Shi's theorem}
Let $q>3$. Then there exists $c>0$ independent of $\Omega$ such that for any
$u_0\in H^{\f 12}_\sigma$ with $\|u_0\|_{H^{\f 12}}\le c$,
the equations \eqref{eq:NSrotating}
admit a unique mild solution $u\in C([0,\infty), H^{\f12}_\sigma)$ satisfying
\ben\label{condition}
\begin{array}{ll}
&u\in C((0,\infty), L^q)\quad\hbox{and}\quad
\lim_{t\rightarrow 0^+}\sup_{0<s<t}s^{\f12-\f3{2q}}\|u(s,\cdot)\|_{L^q}=0,\\
&\na u\in C((0,\infty), L^2)\quad\hbox{and}\quad
\lim_{t\rightarrow 0^+}\sup_{0<s<t}s^{\f14}\|\na u(s,\cdot)\|_{L^2}=0.
\end{array}
\een
Here $H^{\f 12}_\sigma$ denotes the closure of the set
$\{u\in C_c^\infty(\R^3)^3, \dv u=0\}$ in the norm of $\|\cdot\|_{H^{\f 12}}$.
\end{Theorem}

The goal of this paper is to prove the global existence of
(\ref{eq:NSrotating}) for a class of highly oscillating initial
velocity. Thus we need to solve the system \eqref{eq:NSrotating} for
the initial data in a critical functional framework whose regularity
index is negative, for example $\dot B^{-1+\f3p}_{p,q}$ for $p>3$.
However, Cannone's proof \cite{Can1} doesn't work for our case, since
it also relies on the global $L^p-L^q$ estimates for the Stokes
semigroup. Indeed for the Stokes-Coriolis semigroup  $\cG(t)$, one has
\beno
\|\cG(t)u_0\|_{L^p}\le C_{p,\Omega}t^2\|u_0\|_{L^p},\quad \textrm{if} \quad p\neq 2,
\eeno
see Proposition 2.2 in \cite{HS}. Then we can infer from the definition
of Besov space that \beno \|\cG(t)u_0\|_{\dot B^{-1+\f3p}_{p,q}}\le
Ct^2\|u_0\|_{\dot B^{-1+\f3p}_{p,q}}. \eeno
This means that even if the initial data $u_0$ is small in $\dot B^{-1+\f3p}_{p,q}$, the
linear part $\|\cG(t)u_0\|_{\dot B^{-1+\f3p}_{p,q}}$ of the solution
may become large after some time $t_0>0$.

Fortunately, we have the following important observation: if $u\in L^p$
with $\textrm{supp}\,\hat u\in \{\xi: |\xi|\gtrsim \lambda\}$,
then
\beno
\|\cG(t)u\|_{L^p}\le C_{p,\Omega}e^{-t\lambda^2}\|u\|_{L^p}
\eeno
for any $p\in [1,\infty]$ and $t\in [0,\infty]$; while for any $u\in L^2$,
\beno
\|\cG(t)u\|_{L^2}\le \|u\|_{L^2}.
\eeno
This motivates us to introduce the hybrid-Besov spaces $\dot\cB^{\f 12,\f 3p-1}_{2,p}$(see Definition \ref{Def:hybridBes}).
Roughly speaking, if $u\in \dot\cB^{\f 12,\f 3p-1}_{2,p}$, the low frequency part of $u$ belongs to $\dot H^\f12$
and the high frequency part belongs to $\dot B^{-1+\f 3p}_{p,\infty}$. So, $\dot\cB^{\f 12,\f 3p-1}_{2,p}$ is still a
critical space. A remarkable property of $\dot\cB^{\f 12,\f 3p-1}_{2,p}$ is that if $p>3$, then
\beno
\|u_0(x)\|_{\dot\cB^{\f 12,\f 3p-1}_{2,p}}\le C\varepsilon^{1-\f3p},
\eeno
for $u_0(x)=\sin(\f {x_1} \varepsilon)\phi(x), \phi(x)\in \cS(\R^3)$, see Proposition \ref{Prop:oscillate}.
That is, the highly oscillating function is still small in the norm of $\dot\cB^{\f 12,\f 3p-1}_{2,p}$.

Before stating our main results, let us introduce the following space.

\begin{Def}
Let $1\le p\le\infty$, $E_p$ denotes the space of functions such that
$$\|u\|_{E_p}\eqdefa \|u\|_{\widetilde{L}^\infty(\R^+; \dot{\cB}^{\f 12,\f 3p-1}_{2,p})}+
\|u\|_{\widetilde{L}^1(\R^+; \dot{\cB}^{\f 52,\f 3p+1}_{2,p})}<\infty.$$
\end{Def}

Our main results are stated as follows.

\begin{Theorem}\label{TheoremBesov}
Let $p\in [2,4]$. There exists a positive constant $c$ independent of $\Omega$ such that if $\|u_0\|_{\dot\cB^{\f 12,\f 3p-1}_{2,p}}\le
c$, then a unique solution of \eqref{eq:NSrotating} exists in the
small ball of center 0 of $E_p$.
\end{Theorem}

\begin{Remark}
Due to the inclusion map
$$H^{\f12}\subseteq\dot{\cB}^{\f 12,\f 3p-1}_{2,p}\quad\textrm{ for }p\ge 2,$$
it is an improvement of Theorem \ref{Th:Shi's theorem}.
The importance of this improvement is that it allows to construct global solutions of (\ref{eq:NSrotating}) for
a class of  highly oscillating initial velocity $u_0$. For example,
\begin{align*}
u_0(x)=\sin\bigl(\f {x_3} {\varepsilon}\bigr)(-\p_2\phi(x), \p_1\phi(x),0)
\end{align*}
where $\phi\in \cS(\R^3)$ and $\varepsilon>0$ is small enough. This type data is large in the Sobolev norm, however
it is small in the norm of Besov spaces with negative regularity index.
\end{Remark}

\begin{Remark} The inhomogeneous part of the solution has more regularity:
$$
u-\cG(t)u_0\in C(\R^+;\dot B^\f12_{2,\infty}),
$$
which can be proved by following the proof of Proposition \ref{Prop:binear estimate}.
\end{Remark}

If the initial data $u_0\in \dot H^\f12$, we can obtain the following global well-posed result.

\begin{Theorem}\label{Theorem C(H^12)}
Let $p\in [2,4]$. There exists a  positive constant $c$ independent of
$\Omega$ such that if $u_0$ belongs to $\dot{H}^{\f12}$ with
$\|u_0\|_{\dot{\cB}^{\f 12,\f 3p-1}_{2,p}}\le c$, then there exists
a unique global solution of \eqref{eq:NSrotating} in
$C(\R^+,\dot{H}^{\f12})$.
\end{Theorem}
\begin{Remark}  Since we only impose the smallness condition of the initial data in the norm of
$\dot{\cB}^{\f 12,\f 3p-1}_{2,p}$, this
allows to obtain the global well-posedness of (\ref{eq:NSrotating}) for
a class of  highly oscillating initial velocity $u_0$. Moreover, the uniqueness holds in the class
$C(\R^+,\dot{H}^{\f12})$, i.e. it is  unconditional.
\end{Remark}

The structure of this paper is as follows.\vspace{.1cm}

In Section 2, we recall some basic facts about Littlewood-Paley
theory and the functional spaces. In Section 3, we recall some
results concerning the Stoke-Coriolis semigroup's regularizing
effect. Section 4 is devoted to the important bilinear estimates. In
Section 5, we prove Theorem \ref{TheoremBesov} and Theorem
\ref{Theorem C(H^12)}.

\section{Littlewood-paley theory and the function spaces}
First of all, we introduce the Littlewood-Paley decomposition. Choose two
radial functions  $\varphi, \chi \in {\cS}(\mathbb{R}^3)$ supported in
${\cC}=\{\xi\in\mathbb{R}^3,\, \frac{3}{4}\le|\xi|\le\frac{8}{3}\}$,
${\cB}=\{\xi\in\mathbb{R}^3,\, |\xi|\le\frac{4}{3}\}$ respectively such
that
\begin{align*} \sum_{j\in\mathbb{Z}}\varphi(2^{-j}\xi)=1 \quad \textrm{for
all}\,\,\xi\neq 0.
\end{align*}
For $f\in \cS'(\R^3)$, the frequency localization operators $\Delta_j$ and $S_j(j\in\Z)$ are defined by
\begin{align}
\Delta_jf=\varphi(2^{-j}D)f,\quad S_jf=\chi(2^{-j}D)f, \quad D=\frac {\na_x} i.\nonumber
\end{align}
Moreover, we have
$$
S_jf=\sum_{k=-\infty}^{j-1}\Delta_kf\quad\textrm{ in }\quad {\calZ'}(\mathbb{R}^3).
$$
Here we denote the space ${\cal Z'}(\mathbb{R}^3)$ by the dual
space of ${\cal Z}(\mathbb{R}^3)=\{f\in
{\cS}(\mathbb{R}^3);\,D^\alpha \hat{f}(0)=0;
\forall\alpha\in\big(\mathbb{N}\cup 0\big)^3
\,\mbox{multi-index}\}$.

 With our choice of $\varphi$, it is easy to verify that
\begin{align}\label{orth}
\begin{aligned}
&\Delta_j\Delta_kf=0\quad \textrm{if}\quad|j-k|\ge 2\quad
\textrm{and}
\quad \\
&\Delta_j(S_{k-1}f\Delta_k f)=0\quad \textrm{if}\quad|j-k|\ge 5.
\end{aligned}
\end{align}

In the sequel, we will constantly use the Bony's decomposition from \cite{Bony}:
\begin{align}\label{Bonydecom}
fg=T_fg+T_gf+R(f,g), \end{align} with
$$T_fg=\sum_{j\in\mathbb{Z}}S_{j-1}f\Delta_jg,
\quad R(f,g)=\sum_{j\in\mathbb{Z}}\Delta_jf \widetilde{\Delta}_{j}g,
\quad \widetilde{\Delta}_{j}g=\sum_{|j'-j|\le1}\Delta_{j'}g.$$

Let us first recall the definition of general Besov space.

\begin{Def}\label{Def:Bes} Let $s\in\mathbb{R}$, $1\le p,
q\le+\infty$. The homogeneous Besov space $\dot{B}^{s}_{p,q}$ is
defined by
$$\dot{B}^{s}_{p,q}\eqdefa\big\{f\in {\calZ'}(\mathbb{R}^3):\,\|f\|_{\dot{B}^{s}_{p,q}}<+\infty\big\},$$
where
\begin{align*}
\|f\|_{\dot{B}^{s}_{p,q}}\eqdefa \Bigl\|2^{ks}
\|\Delta_kf(t)\|_{L^p}\Bigr\|_{\ell^q}.\end{align*}\end{Def}
If  $p=q=2$,  $\dot{B}^{s}_{2,2}$ is equivalent to the homogeneous Sobelev space $\dot{H}^{s}$.
\vspace{0.1cm}

Now we introduce the hybrid-Besov space  we will work with in this
paper.
\begin{Def}\label{Def:hybridBes} Let $s$,
$\sigma\in\mathbb{R}$, $1\le p\le +\infty$. The hybrid-Besov space
$\dot{\cB}^{s,\sigma}_{2,p}$ is defined by
\begin{align*}
&\dot{\cB}^{s,\sigma}_{2,p}\eqdefa\big\{f\in{\calZ'}(\mathbb{R}^3):
\|f\|_{\dot{\cB}^{s,\sigma}_{2,p}}<+\infty\big\},
\end{align*}
where
$$\|f\|_{\dot{\cB}^{s,\sigma}_{2,p}}\eqdefa\sup_{2^{k}\le \Omega}2^{ks}\|\Delta_k f\|_{L^2}
+\sup_{2^k>\Omega}2^{k\sigma}\|\Delta_k f\|_{L^p}.
$$
\end{Def}

The norm of the space $\widetilde{L}^r_T(\dot{\cB}^{s,\sigma}_{2,p})$ is defined by
$$
\|f\|_{\widetilde{L}^r_T(\dot{\cB}^{s,\sigma}_{2,p})}
\eqdefa\sup_{2^k\le \Omega}2^{ks}\|\Delta_k f\|_{L^r_TL^2}
+\sup_{2^k>\Omega}2^{k\sigma}\|\Delta_k f\|_{L^r_TL^p}.
$$
It is easy to check that
$L^r_T(\dot{\cB}^{s,\sigma}_{2,p})\subseteq \widetilde{L}^r_T(\dot{\cB}^{s,\sigma}_{2,p})$,
where the norm of $L^r_T(\dot{\cB}^{s,\sigma}_{2,p})$ is defined by
$$
\|f\|_{L^r_T(\dot{\cB}^{s,\sigma}_{2,p})}\eqdefa
\big\|\|f(t)\|_{\dot{\cB}^{s,\sigma}_{2,p}}\big\|_{L^r_T}.
$$

The following Berstein's lemma will be repeatedly used throughout this paper.

\begin{Lemma}\cite{Chemin-book}\label{Lem:Bernstein}
Let $1\le p\le q\le+\infty$. Then for any $\beta,\gamma\in(\mathbb{N}\cup\{0\})^3$, there exists a constant $C$
independent of $f$, $j$ such that
\begin{align*} &{\rm supp}\hat f\subseteq
\{|\xi|\le A_02^{j}\}\Rightarrow \|\partial^\gamma f\|_{L^q}\le
C2^{j{|\gamma|}+j n(\frac{1}{p}-\frac{1}{q})}\|f\|_{L^p},
\\
&{\rm supp}\hat f\subseteq \{A_12^{j}\le|\xi|\le
A_22^{j}\}\Rightarrow \|f\|_{L^p}\le
C2^{-j|\gamma|}\sup_{|\beta|=|\gamma|}\|\partial^\beta f\|_{L^p}.
\end{align*}
\end{Lemma}

\begin{Proposition}\label{Prop:oscillate} Let $\phi\in \cS(\R^3)$ and $p>3$.
If $\phi_\varepsilon(x)\eqdefa e^{i\f {x_1} \varepsilon}\phi(x)$, then for any $0<\varepsilon\le \Omega^{-1}$,
\begin{align*}
\|\phi_\varepsilon\|_{\dot\cB^{\f 12,\f 3 p-1}_{2,p}}\le C\varepsilon^{1-\f 3p},
\end{align*}
here $C$ is a constant independent of $\varepsilon$.
\end{Proposition}

\noindent{\bf Proof.}\,
Let $j_0\in\N$ be such that $\Omega\le 2^{j_0}\sim \varepsilon^{-1}$. By Lemma \ref{Lem:Bernstein}, we have
\begin{align*}
\sup_{j\ge j_0}2^{(\f 3p-1)j}\|\Delta_j \phi_\varepsilon\|_{L^p}\le C2^{(\f 3p-1)j_0}\le C\varepsilon^{1-\f 3p}.
\end{align*}
Noting that $e^{i\f {x_1} \varepsilon}=(-i\varepsilon\p_1)^Ne^{i\f {x_1} \varepsilon}$ for any $N\in\N$,
we get by integration by parts that
\begin{align*}
\Delta_j\phi_\varepsilon(x)=(i\varepsilon)^N2^{3j}\int_{\R^3}
e^{i\f {y_1} \varepsilon}\p_{y_1}^N(h(2^j(x-y))\phi(y))dy,\quad h(x)\eqdefa (\cF^{-1}\varphi)(x).
\end{align*}
By Leibnitz formula, we have
\begin{align*}
|\Delta_j\phi_\varepsilon(x)|\le C\varepsilon^N2^{3j}\sum_{k=0}^N2^{kj}\int_{\R^3}
|(\p_{y_1}^kh)(2^j(x-y))||\p_{y_1}^{N-k}\phi(y)|dy,
\end{align*}
from which and Young's inequality, we infer that for $j\ge 0$,
\begin{align*}
\|\Delta_j\phi_\varepsilon\|_{L^q}\le C\varepsilon^N\sum_{k=0}^N2^{kj}2^{3j}\|(\p_{y_1}^kh)(2^jy)\|_{L^1}
\|\p_{y_1}^{N-k}\phi(y)\|_{L^q}\le C\varepsilon^N2^{jN},
\end{align*}
and for $j\le 0$,
\begin{align*}
\|\Delta_j\phi_\varepsilon\|_{L^q}\le C\varepsilon^N\sum_{k=0}^N2^{kj}2^{3j}\|(\p_{y_1}^kh)(2^jy)\|_{L^q}
\|\p_{y_1}^{N-k}\phi(y)\|_{L^1}\le C\varepsilon^N2^{(1-\f 1 q)3j}.
\end{align*}
Thus we have
\begin{align*}
&\sup_{\Omega<2^j<2^{j_0}}2^{(\f 3p-1)j}\|\Delta_j \phi_\varepsilon\|_{L^p}
\le C\varepsilon^N 2^{(N-1+\f 3p)j_0}\le C\varepsilon^{1-\f 3p},\\
&\sup_{2^j\le  \Omega}2^{\f j2}\|\Delta_j \phi_\varepsilon\|_{L^2}\le C\Omega^{\f12}\varepsilon^N
\le C\varepsilon^{N-\f12}.
\end{align*}

Summing up the above estimates yields that
\begin{align*}
\|\phi_\varepsilon\|_{\dot\cB^{\f 12,\f 3 p-1}_{2,p}}\le C\varepsilon^{1-\f 3p}.
\end{align*}
The proof of Proposition \ref{Prop:oscillate} is completed.\endproof

\section{Regularizing effect of the Stokes-Coriolis semigroup}
We consider the linear system
\begin{equation}\label{eq:linsys:fastrot}
\left\{
\begin{array}{ll}
u_t-\nu\Delta u+\Omega e_3\times u+\nabla p=0,\\
\dv u=0, \\
u(0,x)=u_0(x).
\end{array}
\right.
\end{equation}
From Proposition 2.1 in \cite{HS}, we know that
\begin{align}\label{eq:representsol}
\widehat{u}(t,\xi)=\cos\Big(\Omega\f{\xi_3}{|\xi|}t\Big)e^{-\nu|\xi|^2t}I\widehat{u_0}(\xi)+
\sin\Big(\Omega\f{\xi_3}{|\xi|}t\Big)e^{-\nu|\xi|^2t}R(\xi)\widehat{u_0}(\xi),\quad t\ge0,\,\xi\in\R^3,
\end{align}
where $I$ is the identity matrix and
\begin{equation*}\label{}
R(\xi)=\left(
\begin{array}{ccc}
0\, & \,\f{\xi_3}{|\xi|}\, & \,-\f{\xi_2}{|\xi|}\\
-\f{\xi_3}{|\xi|}\, &\, 0 \,& \,\f{\xi_1}{|\xi|} \\
\f{\xi_2}{|\xi|}\, & \,-\f{\xi_1}{|\xi|} \,& \, 0
\end{array}
\right).
\end{equation*}
The Stokes-Coriolis semigroup is explicitly represented by
\begin{align}\label{eq:Green matrix}
\cG(t)f=[\cos(\Omega R_3t)I+\sin(\Omega R_3t)R]e^{\nu t\Delta}f,\quad \textrm{for}\quad t\ge0,\, f\in L^p_\sigma,
\end{align}
where $\widehat{R_3f}(\xi):=\f{\xi_3}{|\xi|}\widehat{f}(\xi)$ for $\xi\neq0$.

We have the following smoothing effect of the Stokes-Coriolis semigroup $\cG$.
\begin{Proposition}\label{Le:Green-Lpest}
Let ${\cC}$ be a ring  centered  at 0 in $\R^3$ . Then there exist  positive  constants
$c$ and $C$ depending only on $\nu$ such that if ${\rm
supp}\,\hat{u}\subset\lambda{\cC}$, then we have\vspace{.2cm}

{\rm (i)}\, for any  $\lambda>0$,
\begin{align}\label{eq:lowGreen}
&\|{\cG}(t)u\|_{L^2}\le Ce^{-c\lambda^2
t}\|u\|_{L^2};
\end{align}

{\rm (ii)}\, if $\lambda\gtrsim \Omega$, then for any $1\le p\le \infty$,
\begin{align}
\|{\cG}(t)u\|_{L^p}\le C e^{-c\lambda^2 t}\|u\|_{L^p}.\label{eq:highGreen}
\end{align}
\end{Proposition}

\noindent{\bf Proof.}\,\,(i)\, Thanks to \eqref{eq:representsol} and Plancherel theorem, we get
\begin{align*}
\|\cG(t)u\|_{L^2}=\|\widehat{\cG}(t,\xi)\hat{u}(\xi)\|_{L^2}\le C
\|e^{-\nu|\xi|^2 t}\hat{u}(\xi)\|_2\le Ce^{-\nu\lambda^2 t}\|u\|_2,
\end{align*}
where  we have used  the support property of $\hat{u}(\xi)$.

(ii)\,Let $\phi\in {\cD}(\mathbb{R}^3\setminus\{0\})$, which equals to 1 near the  ring ${\cC}$.
Set
$$g(t,x)\eqdefa (2\pi)^{-3}\int_{\mathbb{R}^3}e^{ix\cdot\xi}
\phi(\lambda^{-1}\xi)\widehat{\mathcal{G}}(t,\xi)d\xi.
$$
To prove \eqref{eq:highGreen}, it suffices to show
\begin{align}\label{eq:gH1}
\|g(x,t)\|_{L^1}\le Ce^{-c\lambda^2 t}.
\end{align}
Thanks to \eqref{eq:Green matrix}, we infer that
\begin{align}\label{eq:gH2}
\int_{|x|\le\lambda^{-1}}|g(x,t)|dx&\le C
\int_{|x|\le\lambda^{-1}}\int_{\mathbb{R}^3}|\phi(\lambda^{-1}\xi)||\widehat{\mathcal{G}}(t,\xi)|d\xi dx
\le Ce^{-c\lambda^2 t}.
\end{align}
Set $L\eqdefa \frac {x\cdot \nabla_\xi} {i|x|^2}$. Noting that $L(e^{ix\cdot\xi})=e^{ix\cdot\xi}$,
we get by integration by part  that
\begin{align*}
g(x,t)=& \int_{\mathbb{R}^3}L^{N}(e^{ix\cdot\xi})\phi(\lambda^{-1}\xi)
\widehat{\mathcal{G}}(t,\xi)d\xi\nonumber\\=&
\int_{\mathbb{R}^3}e^{ix\cdot\xi}(L^*)^{N}\big(\phi(\lambda^{-1}\xi)
\widehat{\mathcal{G}}(t,\xi)\big)d\xi,
\end{align*}where $N\in\N$ is chosen later. Using the Leibnitz's formula, it is easy to verify that
\beno
&&\big|\pa^{\gamma}(e^{\pm i\Omega \f{\xi_3}{|\xi|}t})\big|\le C|\xi|^{-|\gamma|}(1+ \Omega t)^{|\gamma|},\\
&&\big|\pa^{\gamma}(e^{-\nu|\xi|^2t})\big|\le C|\xi|^{-|\gamma|}e^{-\f \nu 2|\xi|^2t}.
\eeno
Thus we obtain
\begin{align}
&\big|(L^*)^{N}\big(\phi(\lambda^{-1}\xi)\widehat{\mathcal{G}}(t,\xi)\big)\big|\nonumber\\&
\le C|x|^{-N}\sum_{\tiny\left.\begin{array}{cc}|\alpha_1|+|\alpha_2|+|\alpha_3|=|\alpha|\\|\alpha|\le N\end{array}\right.}
\lambda^{-N+\alpha}\big|(\na^{N-\alpha}\phi)(\lambda^{-1}\xi)
\pa^{\alpha_1}(e^{\pm i\Omega \f{\xi_3}{|\xi|}t})\pa^{\alpha_2}(e^{-\nu|\xi|^2t})\pa^{\alpha_3}(I+R(\xi))\big|
\nonumber\\&\le C|\lambda x|^{-N}
\sum_{\tiny\left.\begin{array}{cc}|\alpha_1|+|\alpha_2|+|\alpha_3|=|\alpha|\\|\alpha|\le N\end{array}\right.}
\lambda^{\alpha}|(\na^{N-\alpha}\phi)(\lambda^{-1}\xi)|
|\xi|^{-|\alpha_1|-|\alpha_2|-|\alpha_3|}e^{-\f\nu2|\xi|^2t}(1+ \Omega t)^{|\alpha_1|}.\nonumber
\end{align}
Taking $N=4$, for any $\xi\in \{\xi: A^{-1}\la \le|\xi|\le A\lambda\}$ for some constant $A$ depending on the ring $\cC$
and $\lambda\gtrsim \Omega$,
\begin{align*}
\big|(L^*)^{4}\big(\phi(\lambda^{-1}\xi)\widehat{\mathcal{G}}(t,\xi)\big)\big|\le C|\lambda x|^{-4}e^{-\f \nu 4|\xi|^2t},
\end{align*}
which implies that
\begin{align*}
\int_{|x|\ge\f 1\lambda}|g(x,t)|dx&\le C
e^{-c\lambda^2t}\lambda^3\int_{|x|\ge\f 1\lambda}|\lambda
x|^{-4}dx\le C
e^{-c\lambda^2t},
\end{align*}
which together with \eqref{eq:gH2} gives \eqref{eq:gH1}.
Then the inequality \eqref{eq:highGreen} is proved.\endproof

\vspace{0.1cm}

The following proposition is a direct consequence of Proposition \ref{Le:Green-Lpest}.
\begin{Proposition}\label{prop:semigroup} Let $s,\sigma\in \R$, and $(p,q)\in [1,\infty]$. Then for any
$u\in \dot\cB^{s-\f2q,\sigma-\f2q}_{2,p}$, there holds
\ben\label{eq:hom-est}
\|\cG(t)u\|_{\widetilde{L}^q_T(\dot\cB^{s,\sigma}_{2,p})}\le C\|u\|_{\dot\cB^{s-\f2q,\sigma-\f2q}_{2,p}}.
\een
And for any $f\in \widetilde{L}_T^1\cB^{s,\sigma}_{2,p}$, we have
\ben\label{eq:inhom-est}
\Big\|\int_0^t\cG(t-\tau)f(\tau)d\tau\Big\|_{\widetilde{L}^q_T(\dot\cB^{s+\f2q,\sigma+\f2q}_{2,p})}\le
C\|f(t)\|_{\widetilde{L}^1_T(\dot\cB^{s,\sigma}_{2,p})}.
\een
\end{Proposition}

\no{\bf Proof.}\,\,Here we only prove (\ref{eq:inhom-est}). For any
$2^j\ge \Omega$, we get by Proposition \ref{Le:Green-Lpest} that
\beno \Big\|\Delta_j\int_0^t\cG(t-\tau)f(\tau)d\tau\Big\|_{L^p}\le
C\int_0^te^{-c(t-\tau)2^{2j}}\|\Delta_jf(\tau)\|_{L^p}d\tau, \eeno
from which and Young's inequality, it follows that
\begin{align}\label{eq:High-est}
\Big\|\Delta_j\int_0^t\cG(t-\tau)f(\tau)d\tau\Big\|_{L^q_TL^p}\le&
C\|e^{-ct2^{2j}}\|_{L^q_T}\|\Delta_jf(\tau)\|_{L^1_TL^p}\nonumber\\
\le& C2^{-\f2qj}\|\Delta_jf(\tau)\|_{L^1_TL^p}.
\end{align}
Similarly, we also have
\begin{align}\label{eq:Low-est}
\Big\|\Delta_j\int_0^t\cG(t-\tau)f(\tau)d\tau\Big\|_{L^q_TL^2}\le&
C\|e^{-ct2^{2j}}\|_{L^q_T}\|\Delta_jf(\tau)\|_{L^1_TL^2}\nonumber\\
\le& C2^{-\f2qj}\|\Delta_jf(\tau)\|_{L^1_TL^2}.
\end{align}
Then the inequality (\ref{eq:inhom-est}) follows from (\ref{eq:High-est}) and (\ref{eq:Low-est}).\endproof

\section{Bilinear estimates}

In this section, we study the continuity of the inhomogeneous term in the space $E_{p,T}$ whose norm is defined by
$$\|u\|_{E_{p,T}}\eqdefa \|u\|_{\widetilde{L}^\infty(0,T; \dot{\cB}^{\f 12,\f 3p-1}_{2,p})}+
\|u\|_{\widetilde{L}^1(0,T; \dot{\cB}^{\f 52,\f 3p+1}_{2,p})}.$$
We denote
\beno
B(u,v)\eqdefa\int_0^t\cG(t-\tau)\mathbb{P}\na\cdot(u\otimes v)d\tau,
\eeno
where $\mathbb{P}$ denotes the Helmholtz projection.

\begin{Proposition}\label{Prop:binear estimate}
Let  $p\in [2,4]$. Assume that $u,v\in E_{p,T}$. There exists a
constant $C$ independent of $\Omega$ such that for any $T>0$,
\begin{align}\label{eq:nonhom estimate}
\|B(u,v)\|_{E_{p,T}}\le C\|u\|_{E_{p,T}}\|v\|_{E_{p,T}}.
\end{align}

\end{Proposition}
\noindent{\bf Proof.}\,\,
Thanks to  Proposition \ref{prop:semigroup}, it suffices to show that
\begin{align}\label{eq:bilinear estimate}
\|uv\|_{\widetilde{L}_T^1\dot\cB^{\f 32,\f 3p}_{2,p}}\le C\|u\|_{E_{p,T}}\|v\|_{E_{p,T}}.
\end{align}

From Bony's decomposition \eqref{Bonydecom} and \eqref{orth}, we have
\begin{align*}\label{}
\Delta_j(u v)=&\sum_{|k-j|\le 4}\Delta_j(S_{k-1}u\Delta_kv)
+\sum_{|k-j|\le 4}\Delta_j(S_{k-1}v \Delta_ku)
\nonumber\\&+\sum_{k\ge
j-2}\Delta_j(\Delta_ku\widetilde{\Delta}_kv)\nonumber\\\eqdefa &
I_j+II_j+III_j.
\end{align*}

Set $J_j\eqdefa \big\{(k',k);\,|k-j|\le 4, k'\le k-2\big\}$,
then for $2^j>\Omega$,
\begin{align*}
\|I_j\|_{L_T^1 L^p}&\le\sum_{J_j}\|\Delta_j(\Delta_{k'}u\Delta_{k}v)\|_{L^1_TL^p}\nonumber\\
&\le\Big(\sum_{J_{j,\ell \ell}}+\sum_{J_{,\ell h}}+\sum_{J_{j,hh}}\Big)
\|\Delta_j(\Delta_{k'}u\Delta_{k}v)\|_{L_T^1 L^p}\\
&\eqdefa I_{j,1}+I_{j,2}+I_{j,3},
\end{align*}
where
\begin{align*}
&J_{j,\ell \ell}=\{(k',k)\in J_j,\,\, 2^{k'}\le \Omega, 2^k\le \Omega\},\\&
J_{j,\ell h}=\{(k',k)\in J_j,\,\, 2^{k'}\le \Omega, 2^k>\Omega\},\\
&J_{j,hh}=\{(k',k)\in J_j,\,\, 2^{k'}>\Omega, 2^k>\Omega\}.\end{align*}
We get by using Lemma \ref{Lem:Bernstein}  that
\begin{align*}
I_{j,1}&\le C\sum_{(k',k)\in J_{j,\ell\ell}}\!\!
\|\Delta_{k'}u\|_{L_T^{\infty}L^\infty}2^{k(\f 32-\f 3p)}\|\Delta_{k}v\|_{L_T^{1}L^2}\nonumber\\&\le C
\sum_{(k',k)\in J_{j,\ell \ell}}2^{\f {k'}2}\|\Delta_{k'}u\|_{L_T^{\infty}L^2}
2^{k'}\|\Delta_{k}v\|_{L_T^{1}L^2}
2^{k(\f 32-\f 3p)}\\
&\le C \|u\|_{\widetilde{L}_T^{\infty}\dot\cB^{\f 12,\f 3p-1}_{2,p}}
\|v\|_{\widetilde{L}_T^{1}\dot\cB^{\f 52,\f 3p+1}_{2,p}}\sum_{(k',k)\in J_{j,\ell \ell}}2^{(k'-k)}2^{-\f3pk}\\
&\le C2^{-\f 3{p}j}\|u\|_{\widetilde{L}_T^{\infty}\dot\cB^{\f 12,\f 3p-1}_{2,p}}
\|v\|_{\widetilde{L}_T^{1}\dot\cB^{\f 52,\f 3p+1}_{2,p}},
\end{align*}
where we used in the last inequality the fact that
\beno
\sum_{(k',k)\in J_{j,\ell \ell}}2^{(k'-k)}2^{-\f3pk}\le \sum_{k'\le k-2}2^{(k'-k)}\sum_{|k-j|\le 4}2^{-\f3pk}
\le C2^{-\f3pj},
\eeno
with $C$ independent of $j$. Similarly, we have
\begin{align*}
I_{j,2}&\le \sum_{(k',k)\in J_{j,\ell h}}\|\Delta_{k'}u\|_{L_T^{\infty}L^\infty}
\|\Delta_{k}v\|_{L_T^{1}L^p}\\
&\le C\sum_{(k',k)\in J_{j,\ell h}}2^{\f {k'}2}\|\Delta_{k'}u\|_{L_T^{\infty}L^2}
2^{k'}\|\Delta_{k}v\|_{L_T^{1}L^p}\nonumber\\&
\le C2^{-\f 3{p}j}\|u\|_{\widetilde{L}_T^{\infty}\dot\cB^{\f 12,\f 3p-1}_{2,p}}
\|v\|_{\widetilde{L}_T^{1}\dot\cB^{\f 52,\f 3p+1}_{2,p}},
\end{align*}
and
\begin{align*}
I_{j,3}&\le \sum_{(k',k)\in J_{j,hh}}\|\Delta_{k'}u\|_{L_T^{\infty}L^\infty}
\|\Delta_{k}v\|_{L_T^{1}L^p}\\
&\le C\sum_{(k',k)\in J_{j,hh}}2^{k'(\f 3p-1)}\|\Delta_{k'}u\|_{L_T^{\infty}L^p}
2^{k'}
\|\Delta_{k}v\|_{L_T^{1}L^p}\nonumber\\&
\le C2^{-\f 3{p}j}
\|u\|_{\widetilde{L}_T^{\infty}\dot\cB^{\f 12,\f 3p-1}_{2,p}}
\|v\|_{\widetilde{L}_T^{1}\dot\cB^{\f 52,\f 3p+1}_{2,p}}.
\end{align*}

On the other hand, for $2^j\le \Omega$, we have
\begin{align*}
\|I_j\|_{L_T^1L^2}&\le
\sum_{{J_j}}\|\Delta_j(\Delta_{k'}u\Delta_{k}v)\|_{L_T^1L^2}
\nonumber\\&\le \Big(\sum_{{J}_{j,\ell \ell}}+\sum_{{J}_{j,\ell h}}
+\sum_{{J}_{j,hh}}\Big)
\|\Delta_j(\Delta_{k'}u\Delta_{k}v)\|_{L_T^{1}L^2}\\
&\eqdefa I_{j,4}+I_{j,5}+I_{j,6}.
\end{align*}
We get by using Lemma \ref{Lem:Bernstein} that
\begin{align*}
I_{j,4}&\le C\sum_{(k,k')\in J_{j,\ell \ell}}2^{\f {k'}2}\|\Delta_{k'}u\|_{L^{\infty}_TL^2}2^{k'}
\|\Delta_{k}v\|_{L^{1}_TL^2}\\&\le C2^{-\f {3j}{2}}\|u\|_{\widetilde{L}_T^{\infty}\dot\cB^{\f 12,\f 3p-1}_{2,p}}
\|v\|_{\widetilde{L}_T^{1}\dot\cB^{\f 52,\f 3p+1}_{2,p}},
\end{align*}
and noting  $p\le 4,$
\begin{align*}
I_{j,5}&\le
C\sum_{(k,k')\in {J}_{j,\ell h}}\|\Delta_{k'}u\|_{L^{\infty}_TL^{\f {2p} {p-2}}}
\|\Delta_{k}v\|_{L^{1}_TL^p}\\ &\le
C\sum_{(k,k')\in {J}_{j,\ell h}}2^{\f {k'}2}\|\Delta_{k'}u\|_{L^{\infty}_TL^2}2^{k'(\f 3p-\f 12)}
\|\Delta_{k}v\|_{L^{1}_TL^p}\\&
\le C 2^{-\f {3j}{2}}
\|u\|_{\widetilde{L}_T^{\infty}\dot\cB^{\f 12,\f 3p-1}_{2,p}}
\|v\|_{\widetilde{L}_T^{1}\dot\cB^{\f 52,\f 3p+1}_{2,p}},
\end{align*}
and
\begin{align*}
I_{j,6}&\le C
\sum_{(k,k')\in J_{j,hh}}\|\Delta_{k'}u\|_{L^{\infty}_TL^{\f {2p} {p-2}}}
\|\Delta_{k}v\|_{L^{1}_TL^p}\nonumber\\
&\le C\sum_{(k,k')\in J_{j,hh}}2^{k'(\f 3p-1)}\|\Delta_{k'}u\|_{L^{\infty}_TL^p}2^{k'(\f {3}p-\f 12)}
\|\Delta_{k}v\|_{L^{1}_TL^p}\nonumber\\&\le C2^{-\f {3j}{2}}
\|u\|_{\widetilde{L}_T^{\infty}\dot\cB^{\f 12,\f 3p-1}_{2,p}}
\|v\|_{\widetilde{L}_T^{1}\dot\cB^{\f 52,\f 3p+1}_{2,p}}.
\end{align*}

Summing up the above estimates $I_{j,1}-I_{j,6}$ yields that
\begin{align}\label{eq:paraproduct of convect term1}
&\sup_{2^j>1}2^{j\f 3p}\|I_j\|_{L^1_TL^p}+ \sup_{2^j\le1}2^{\f
{3j}2}\|I_j\|_{L^1_TL^2}\le C\|u\|_{E_{p,T}}\|v\|_{E_{p,T}}.
\end{align}
By the same procedure as to \eqref{eq:paraproduct of convect term1}, we have
\begin{align}\label{eq:paraproduct of convect term2}
&\sup_{2^j>1}2^{j\f 3p}\|II_j\|_{L^{1}_TL^p}+
\sup_{2^j\le1}2^{\f {3j}2}\|II_j\|_{L^1_TL^2}\le C\|u\|_{E_{p,T}}\|v\|_{E_{p,T}}.
\end{align}

Set ${K_j}\eqdefa \{(k,k');\,k\ge j-3, |k'-k|\le1\}$.
Then we have
\begin{align*}
III_j&=\Bigl(\sum_{K_{j,\ell\ell}}+\sum_{K_{j,\ell h}}+\sum_{K_{j,h\ell}}+\sum_{K_{j,hh}}\Bigr)
\Delta_j(\Delta_{k}u{\Delta}_{k'}v)\\
&\eqdefa III_{j,1}+III_{j,2}+III_{j,3}+III_{j,4},
\end{align*}
where
\begin{align*}
&K_{j,\ell\ell}=\{(k,k')\in K_j,\, 2^k\le \Omega, 2^{k'}\le \Omega\},\\&
K_{j,\ell m}=\{(k,k')\in K_j,\, 2^k\le \Omega, 2^{k'}>\Omega\},\\&
K_{j,hm}=\{(k,k')\in K_j,\, 2^k>\Omega, 2^{k'}\le \Omega\},\\&
K_{j,hh}=\{(k,k')\in K_j,\, 2^k>\Omega, 2^{k'}>\Omega\}.\end{align*}
We get by Lemma  \ref{Lem:Bernstein}   that
\begin{align*}
&\|III_{j,1}\|_{L^1_TL^p}\le C2^{3j(1-\f 1{p})}\!\sum_{(k,k')\in K_{j,\ell\ell}}
\|\Delta_ku\Delta_{k'}v\|_{L^{1}_TL^1}\\
&\le C2^{3j(1-\f 1{p})}\!\sum_{(k,k')\in K_{j,\ell\ell}}\!2^{\f k2}
\|\Delta_ku\|_{L^{\infty}_TL^2}
2^{-\f k2}2^{k'\f 52}\|\Delta_{k'}v\|_{L^{1}_TL^2}
2^{-k'\f 52}\nonumber\\
&\le C2^{3j(1-\f 1{p})}\|u\|_{\widetilde{L}_T^{\infty}\dot\cB^{\f 12,\f 3p-1}_{2,p}}
\|v\|_{\widetilde{L}_T^{1}\dot\cB^{\f 52,\f 3p+1}_{2,p}}\!\sum_{(k,k')\in K_{j,\ell\ell}}
2^{-\f k2-\f52k'}\nonumber\\
&\le C2^{-\f 3{p}j}\|u\|_{\widetilde{L}_T^{\infty}\dot\cB^{\f 12,\f 3p-1}_{2,p}}
\|v\|_{\widetilde{L}_T^{1}\dot\cB^{\f 52,\f 3p+1}_{2,p}}\sum_{k\ge j-3}2^{-3(k-j)}\nonumber\\
&\le C2^{-\f 3{p}j}\|u\|_{E_{p,T}}\|v\|_{E_{p,T}},
\end{align*}
and
\begin{align*}
\|III_{j,1}\|_{L^{1}_TL^2}&\le C2^{\f {3j}{2}}\!\sum_{(k,k')\in K_{j,\ell\ell}}
\|\Delta_ku\Delta_{k'}v\|_{L^{1}_TL^1}\\
&\le C2^{-\f {3j}{2}}\|u\|_{E_{p,T}}\|v\|_{E_{p,T}}.
\end{align*}
Similarly,  we obtain
\begin{align*}
\|&III_{j,2}+III_{j,3}\|_{L^1_TL^p}\le C2^{\f {3j}{2}}\!\sum_{(k,k')\in K_{j,\ell h}\cup K_{j,h \ell}}
\|\Delta_ku\Delta_{k'}v\|_{L^{1}_TL^\f{2p}{2+p}}\\
&\le C2^{\f {3j}{2}}\Big(\sum_{K_{j,\ell h}}\|\Delta_ku\|_{L^{\infty}_TL^2}
\|\Delta_{k'}v\|_{L^{1}_TL^p}+
\sum_{K_{j,h\ell}}\|\Delta_ku\|_{L^{1}_TL^p}
\|\Delta_{k'}v\|_{L^{\infty}_TL^2}\Big)
\nonumber\\&\le C2^{-\f 3{p}j}\|u\|_{E_{,T}}\|v\|_{E_{p,T}},
\end{align*}
and
\begin{align*}
\|III_{j,2}+III_{j,3}\|_{L^{1}_TL^2}&\le C2^{\f 3{p}j}\!\sum_{(k,k')\in K_{j,\ell h}\cup K_{j,h\ell}}
\|\Delta_ku\Delta_{k'}v\|_{L^{1}_T{L^\f{2p}{2+p}}}\\
&\le C2^{-\f {3j}{2}}\|u\|_{E_{p,T}}\|v\|_{E_{p,T}}.
\end{align*}
Finally, due to  $2\le p\le 4$, we have
\begin{align*}
\|III_{j,4}\|_{L^1_TL^p}&\le C2^{\f 3{p}j}\!\sum_{(k,k')\in K_{j,hh}}
\|\Delta_ku\Delta_{k'}v\|_{L^{1}_TL^\f{p}{2}}\\
&\le C2^{\f 3{p}j}\!\sum_{(k,k')\in K_{j,hh}}\!\|\Delta_ku\|_{L^{\infty}_TL^p}
\|\Delta_{k'}v\|_{L^{1}_TL^p}
\nonumber\\&\le C2^{-\f 3{p}j}\|u\|_{E_{p,T}}\|v\|_{E_{p,T}},\\
\|III_{j,4}\|_{L^{1}_TL^2}&\le C2^{3j(\f{2}p-\f 12)}\!\sum_{(k,k')\in K_{j,hh}}
\|\Delta_ku\Delta_{k'}v\|_{L^{1}_T{L^\f{p}{2}}}\\
&\le C2^{-\f {3j}2}\|u\|_{E_{p,T}}\|v\|_{E_{p,T}}.
\end{align*}

Summing up the estimates of $III_{j,1}-III_{j,4}$, we obtain
\begin{align}\label{eq:paraproduct of convect term}
\sup_{2^j>1}2^{\f 3pj}\|III_j\|_{L^{1}_TL^p}+
\sup_{2^j\le1}2^{\f {3j}2}\|III_j\|_{L^{1}_TL^2}\le C\|u\|_{E_{p,T}}\|v\|_{E_{p,T}}.
\end{align}

Then the inequality \eqref{eq:bilinear estimate}  can be deduced from
\eqref{eq:paraproduct of convect term1}--\eqref{eq:paraproduct of convect term}.\endproof

In order to  prove the uniqueness of the solution in $C(\R^+;\dot H^\f12)$, we establish
the following new bilinear estimate in the weight time-space Besov space introduced in \cite{CMZ1, CMZ2}.

\begin{Proposition}\label{Prop:uniquebilinearestimate} Assume that $u,v\in L^\infty_T(\dot B^{\f12}_{2,\infty})$.
Then for any $T>0$, we have
\begin{align*}
\|B(u,v)\|_{L^\infty_T\dot{B}^{\f12}_{2,\infty}}\le C
\|u\|_{L^\infty_T\dot{B}^{\f12}_{2,\infty}}
\big\|\omega_{j,T}2^{\f j2}\|\Delta_jv\|_{L^\infty_TL^2}\big\|_{\ell^\infty},
\end{align*}
where
$$\omega_{j,T}\eqdefa \sup_{k\ge j}e_{k,T}2^{\f12(j-k)},\quad  e_{j,T}\eqdefa 1-e^{-c2^{2j}T}.$$
\end{Proposition}

\begin{Remark} Compared with $e_{j,T}$, the weight $\omega_{j,T}$ satisfies (\ref{omega-properties}),
which is important to the following estimates. On the other hand, due to the fact $\lim_{T\rightarrow0}\omega_{j,T}=0$,
it can be proved that if $u\in C([0,T];\dot H^\f12)$,
then for any $\varepsilon>0$, one has
\beno
\big\|\omega_{j,T}2^{\f j2}\|\Delta_jv\|_{L^\infty_TL^2}\big\|_{\ell^\infty}<\varepsilon
\quad \textrm{if T small enough}.
\eeno
This point is important in the proof of the uniqueness.
\end{Remark}

\noindent{\bf Proof}.\,\,Firstly, we list some useful properties of $\omega_{j,T}$:
\ben\label{omega-properties}
\begin{array}{l}
e_{j,T}\le \omega_{j,T}\qquad \textrm{for any}\,\,j\in\Z,\\
\omega_{j,T}\le 2^{\f12(j-j')}\omega_{j',T}\quad\mbox{if}\quad j'\le j
\quad \textrm{and}\quad  \omega_{j,T}\le 2\omega_{j',T} \quad\mbox{if}\quad j\le j'.
\end{array}
\een
We get by Proposition \ref{Le:Green-Lpest} that
\begin{align}\label{eq:uniq:uv}
\|B(u,v)\|_{\dot{B}^{\f12}_{2,\infty}}
&\le\sup_{j\in\Z}2^{\f j2}\int_0^t\|\cG(t-\tau)\Delta_j\mathbb{P}\na\cdot(u\otimes v)\|_{L^2}d\tau\nonumber\\&
\le\sup_{j\in\Z}2^{\f {3j}2}\|e^{-c2^{2j}t}\|_{L^1_T}\|\Delta_j(u\otimes v)\|_{L^\infty_TL^2}\nonumber\\&
\le C\sup_{j\in\Z}2^{-\f j2}e_{j,T}\|\Delta_j(u v)\|_{L^\infty_TL^2}.
\end{align}
We use Bony's decomposition to estimate $\|\Delta_j(u v)\|_{L^\infty_TL^2}$.
Thanks to (\ref{omega-properties}), we have
\begin{align}\label{eq:uniq:Tuv}
\sum_{|k-j|\le4}\|\Delta_j(S_{k-1}u \Delta_kv)\|_{L^\infty_TL^2}
&\le C\|u\|_{L^\infty_T\dot B^{\f12}_{2,\infty}}\sum_{|k-j|\le4}2^k\|\Delta_kv\|_{L^\infty_TL^2}\nonumber\\
&\le C\omega_{j,T}^{-1}2^{\f j2}\|u\|_{L^\infty_T\dot B^{\f12}_{2,\infty}}
\big\|\omega_{k,T}2^{\f {k}2}\|\Delta_{k}v\|_{L^\infty_TL^2}\big\|_{\ell^\infty},
\end{align}
and by (\ref{omega-properties}) again,
\begin{align*}
\|S_{k-1}v\|_{L^\infty}&\le\sum_{k'\le k-2}\|\Delta_{k'}v\|_{L^2}2^{\f32k'}\le
\big\|\omega_{k',T}2^{\f {k'}2}\|\Delta_{k'}v\|_{L^\infty_TL^2}\big\|_{\ell^\infty}
\sum_{k'\le k-2}2^{k'}\omega_{k',T}^{-1}\\
&\le 2^k\omega_{k,T}^{-1}\big\|\omega_{k',T}2^{\f {k'}2}\|\Delta_{k'}v\|_{L^\infty_TL^2}\big\|_{\ell^\infty},
\end{align*}
which implies that
\begin{align}\label{eq:uniq:uTv}
\sum_{|k-j|\le4}&\|\Delta_j(S_{k-1}v \Delta_ku)\|_{L^\infty_TL^2}
\le 2^{\f k2}\omega_{k,T}^{-1}\|u\|_{L^\infty_T\dot B^{\f12}_{2,\infty}}
\big\|\omega_{k',T}2^{\f {k'}2}\|\Delta_{k'}v\|_{L^\infty_TL^2}\big\|_{\ell^\infty},
\end{align}
and for the remainder term,
\begin{align}\label{eq:uniq:Ruv}
\sum_{k\ge j-2}\|\Delta_j(\Delta_{k}u \widetilde{\Delta}_kv)\|_{L^\infty_TL^2}&\le
\sum_{k\ge j-2}2^{\f32j}\|\Delta_j(\Delta_{k}u \widetilde{\Delta}_kv)\|_{L^\infty_TL^1}\nonumber\\&
\le C\sum_{k\ge j-2}2^{\f 32j}\|\Delta_{k}u\|_{L^\infty_TL^2}\|\widetilde{\Delta}_{k}v\|_{L^\infty_TL^2}
\nonumber\\&\le C\omega_{j,T}^{-1}2^{\f j2}\|u\|_{L^\infty_T\dot B^{\f12}_{2,\infty}}
\big\|\omega_{k,T}2^{\f {k}2}\|\Delta_{k}v\|_{L^\infty_TL^2}\big\|_{\ell^\infty}.
\end{align}
Substituting \eqref{eq:uniq:Tuv}-\eqref{eq:uniq:Ruv} into \eqref{eq:uniq:uv}
concludes the proof. \endproof

\section{Proof of Theorem \ref{TheoremBesov} and Theorem \ref{Theorem C(H^12)}}

The proof of Theorem \ref{TheoremBesov} is based on the following classical lemma.
\begin{Lemma}\cite{Can-book}\label{Lem:fixpoint}
Let X be an abstract Banach space and $B:X\times X\rightarrow X$ a
bilinear operator, $\|\cdot\|$ being the $X$-norm, such that for
any $x_1\in X$ and $x_2\in X$, we have
$$
\|B(x_1,x_2)\| \leq \eta \|x_1\|\|x_2\|,
$$
then for any $y\in X$ such that
$$
4\eta\|y\|<1,
$$
the equation
$$
x=y+B(x,x)
$$
has a solution $x$ in $X$. Moreover, this solution $x$ is the only
one such that
$$
\|x\|\leq \frac {1-\sqrt{1-4\eta\|y\|}} {2\eta}.
$$
\end{Lemma}

\noindent{\bf Proof of Theorem \ref{TheoremBesov}.}\,Using the Stokes-Coriolis semigroup,
we rewrite the system  \eqref{eq:NSrotating} as the following integral form
\begin{align}\label{eq:integral form}
u(x,t)=\cG(t)u_0-\int_0^t \cG(t-\tau)\mathbb{P}\na\cdot (u\otimes u) d\tau\eqdefa \cG(t)u_0+B(u,u).
\end{align}
Thanks to Proposition \ref{prop:semigroup}, we have
\beno
\|\cG(t)u_0\|_{E_p}\le C\|u_0\|_{\dot{\cB}^{\f 12,\f 3p-1}_{2,p}}\le Cc.
\eeno
Obviously, $B(u,v)$ is bilinear, and we get by Proposition \ref{Prop:binear estimate} that
\beno
\|B(u,v)\|_{E_p}\le C\|u\|_{E_p}\|v\|_{E_p}.
\eeno
Taking $c$ such that $4C^2c<\f34$, Lemma \ref{Lem:fixpoint} ensures that the equation
\beno
u=\cG(t)u_0+B(u,u)
\eeno
has a unique solution in the ball $\{u\in E_p: \|u\|_{E_p}\le \f 1 {4C}\}$.\endproof

Now we turn to prove Theorem \ref{Theorem C(H^12)}.\vspace{0.2cm}

\noindent{\bf Proof of Theorem \ref{Theorem C(H^12)}.}\,\,
We introduce Banach space $F_p$ whose norm is defined by
\beno
\|u\|_{F_p}\eqdefa \|u\|_{\widetilde{L}^\infty(\R^+;\dot H^\f12)}+\|u\|_{E_p}.
\eeno

{\bf Step 1.}\, Existence in $F_p$

We define the map
\beno
\mathcal{T}u\eqdefa\cG(t)u_0+B(u,u).
\eeno
Next we prove that if $c$ is small enough, the map $\mathcal{T}$ has a unique fixed point in the ball
$$B_A\eqdefa\big\{u\in F_p: \|u\|_{E_p}\le Ac, \|u\|_{F_p}\le A\|u_0\|_{\dot H^\f12}\big\},$$
for some $A>0$ to be determined later.
From Proposition \ref{prop:semigroup} and Proposition \ref{Prop:binear estimate}, we infer that
\ben\label{eq:E_p-est}
\|\mathcal{T}u\|_{E_p}\le C\|u_0\|_{\cB^{\f12,\f3p-1}_{2,p}}+C\|u\|_{E_p}^2.
\een

On the other hand,  we get by Proposition \ref{Le:Green-Lpest} that
\begin{align}\label{eq:regular:uv}
\|B(u,u)\|_{\widetilde{L}^\infty(\R^+;\dot{H}^{\f12})}&
\le\Big\|\int_0^t\cG(t-\tau)\mathbb{P}\na\cdot(u\otimes u)(\tau)d\tau\Big\|_{\widetilde{L}^\infty(\R^+;\dot{H}^{\f12})}\nonumber\\
&\le C\Big(\sum_{j\in\Z}2^j\Big(\sup_{t\in \R^+}
\int_0^t\|\cG(t-\tau)\Delta_j\mathbb{P}\na\cdot(u\otimes u)(\tau)\|_{L^2}d\tau\Big)^2\Big)^{\f12}
\nonumber\\&\le C
\Big\|2^{\f32j}\sup_{t\in\R^+}\int_0^te^{-c2^{2j}t}\|\Delta_j(u\otimes u)\|_{L^2}d\tau\Big\|_{\ell^2}.
\end{align}
In the following, we denote $\{c_j\}_{j\in \Z}$ by a sequence in $\ell^2$
with the norm $\|\{c_j\}\|_{\ell^2(\Z)}\le 1$. We get by Lemma \ref{Lem:Bernstein} that
\begin{align}\label{eq:regular:Tuv}
\sup_{t\in\R^+}\int_0^te^{-c2^{2j}t}\|\Delta_j(T_uu)\|_{L^2}d\tau
&\le\|e^{-c2^{2j}t}\|_{L^1(\R^+)}\sum_{|k-j|\le4}\|\Delta_j(S_{k-1}u\Delta_k u)\|_{L^\infty(\R^+;L^2)}\nonumber\\
&\le C2^{-2j}\|S_{k-1}u\|_{L^\infty(\R^+;L^\infty)}\sum_{|k-j|\le4}
\|\Delta_k u\|_{L^\infty(\R^+;L^2)}\nonumber\\
&\le C\|u\|_{\widetilde{L}^\infty(\R^+;\dot{\cB}^{\f 12,\f 3p-1}_{2,p})}2^k2^{-2j}\sum_{|k-j|\le4}
\|\Delta_k u\|_{L^\infty(\R^+;L^2)}\nonumber\\&\le C2^{-\f32j}\|u\|_{E_p}\sum_{|k-j|\le4}
2^{\f{(k-j)}2}2^{\f k2}\|\Delta_k u\|_{L^\infty(\R^+;L^2)}\nonumber\\
&\le C2^{-\f32j} c_j\|u\|_{E_p}\|u\|_{\widetilde{L}^\infty(\R^+;\dot H^{\f12})}.
\end{align}
The remainder term of $uv$ is estimated by
\begin{align}\label{eq:regular:Ruv}
\sup_{t\in\R^+}\int_0^te^{-c2^{2j}t}\|\Delta_jR(u,u)\|_{L^2}d\tau&\le
\|e^{-c2^{2j}t}\|_{L^\infty(\R^+)}\sum_{k\ge j-2}\|\Delta_j(\Delta_{k}u\widetilde{\Delta}_k u)\|_{L^1(\R^+;L^2)}
\nonumber\\&\le C\sum_{k\ge j-2}\|\widetilde{\Delta}_k u\|_{L^1(\R^+;L^\infty)}
\|\Delta_{k}u\|_{L^\infty(\R^+;L^2)}\nonumber\\
&\le C\|u\|_{\widetilde{L}^1\dot{\cB}^{\f 52,\f 3p+1}_{2,p}}\sum_{k\ge j-2}2^{-k}\|\Delta_{k}u\|_{L^\infty(\R^+;L^2)}
\nonumber\\&\le C\|u\|_{E_p}\sum_{k\ge j-2}2^{-\f32k}2^{\f12k}\|\Delta_{k}u\|_{L^\infty(\R^+;L^2)}\nonumber\\
&\le C2^{-\f32j} c_j\|u\|_{E_p}\|u\|_{\widetilde{L}^\infty(\R^+;\dot H^{\f12})}.
\end{align}
Combining \eqref{eq:regular:Tuv}-\eqref{eq:regular:Ruv} with \eqref{eq:regular:uv} yields that
\begin{align*}
\|B(u,u)\|_{\widetilde{L}^\infty(\R^+;\dot H^{\f12})}\le C\|u\|_{E_p}\|u\|_{\widetilde{L}^\infty(\R^+;\dot H^{\f12})}.
\end{align*}
It is easy to verify that
\begin{align*}
\|\cG(t)u_0\|_{\widetilde{L}^\infty_T\dot{H}^{\f12}} \le C\|u_0\|_{\dot H^{\f12}}.
\end{align*}
Consequently by (\ref{eq:E_p-est}) and $\|u_0\|_{\dot\cB^{\f12,\f3p-1}_{2,p}}\le C\|u_0\|_{\dot H^{\f12}}$
(from Lemma \ref{Lem:Bernstein} and the definition of Besov space), we obtain
\ben\label{eq:F_p-est}
\|\mathcal{T}u\|_{F_p}\le C\|u_0\|_{\dot H^\f12}+C\|u\|_{E_p}\|u\|_{F_p}.
\een
Taking $A=2C$ and $c>0$ such that $2C^2c\le \f12$, it follows from (\ref{eq:E_p-est})
and (\ref{eq:F_p-est}) that the map $\mathcal{T}$ is a map from $B_A$ to $B_A$. Similarly,
it can be proved that $\mathcal{T}$ is also a contraction in $B_A$. Thus, Banach fixed point theorem
ensures that the map $\mathcal{T}$ has a unique fixed point in $B_A$.

{\bf Step 2.}\, Uniqueness in $C(\R^+;\dot H^{\f12})$\,\,

Let $u_1$ and $u_2$ be two solutions of (\ref{eq:NSrotating}) in $C(\R^+;\dot H^{\f12})$
with the same initial data $u_0$. We consider
\begin{align*}
u_1-u_2=&B\big(u_1-\cG(t)u_0, u_1-u_2\big)+B\big(\cG(t)u_0, u_1-u_2\big)\nonumber\\&
+B\big(u_1-u_2, u_2-\cG(t)u_0\big)+B\big(u_1-u_2, \cG(t)u_0\big).
\end{align*}
Then we get by Proposition \ref{Prop:uniquebilinearestimate} that
\begin{align}\label{eq:diff-equ}
\sup_{t\in [0,T]}\|(u_1-u_2)(t)\|_{\dot B^{\f12}_{2,\infty}}
\le& C\sup_{t\in [0,T]}\|(u_1-u_2)(t)\|_{\dot B^{\f12}_{2,\infty}}
\Big(\big\|\omega_{j,T}2^{\f j2}\|\Delta_j u_0\|_2\big\|_{\ell^\infty}
\nonumber\\&+\sup_{t\in [0,T]}\|u_1(t)-\cG(t)u_0\|_{\dot H^{\f12}}
+\sup_{t\in [0,T]}\|u_2(t)-\cG(t)u_0\|_{\dot H^{\f12}}\Big),
\end{align}
here we used the fact $\omega_{j,T}\le 1$ so that
\beno
\big\|\omega_{j,T}2^{\f j2}\|\Delta_ju\|_{L^\infty_TL^2}\big\|_{\ell^\infty}
\le \sup_{t\in [0,T]}\|u(t)\|_{\dot H^{\f12}}.
\eeno

Noticing that $\omega_{j,0}=0$ and $u_0\in\dot H^{\f12}$, we have
\begin{align*}\label{}
\big\|\omega_{j,T}2^{\f j2}\|\Delta_j u_0\|_2\big\|_{\ell^\infty}\le\f1{3C},
\end{align*}
for $T$ small enough. On the other hand, since $u_1,u_2\in
C(\R^+;\dot H^{\f12})$, we also have \beno \sup_{t\in
[0,T]}\|u_1-\cG(t)u_0\|_{\dot H^{\f12}} +\sup_{t\in
[0,T]}\|u_2-\cG(t)u_0\|_{\dot H^{\f12}}\le \f 1 {3C}, \eeno
for $T$ small enough. Then (\ref{eq:diff-equ}) ensures that $u_1(t)=u_2(t)$ for $T$ small enough.
Then by the standard continuity argument, we can conclude that $u_1=u_2$ on $[0,\infty)$.
\endproof

\bigskip

\noindent {\bf Acknowledgments.}
 Q. Chen is supported by  NSF of China under Grants 10701012 and 10931001.
C. Miao is supported by  NSF of China under Grant 10725102. Z. Zhang is
supported by NSF of China under Grant 10990013.


\begin{thebibliography}{50}
\bibitem{BMN1} A. Babin, A. Mahalov and B. Nicolaenko, {\it Regularity and integrability for the 3D Euler
and Navier-Stokes equations for uniformly rotating fluids}. Asympt. Anal., 15(1997), 103-150.

\bibitem{BMN2} A. Babin, A. Mahalov and B. Nicolaenko, {\it Global regularity of 3D rotating Navier-Stokes
equations for resonant domains}. Indiana Univ. Math. J.,  48(1999), 1133-1176.

\bibitem{Bony} J.-M. Bony, {\it Calcul symbolique et propagation
des singulariti\'{e}s pour les \'{e}quations aux d\'{e}riv\'{e}es
partielles non lin\'{e}aires}. Ann. de l'Ecole Norm. Sup., 14(1981), 209-246.

\bibitem{Can-book} M. Cannone, {\it Ondellettes, paraproduits et Navier-Stokes}, Paris: Diderot Editeur, 1995.

\bibitem{Can1} M. Cannone,
{\it A generalization of a theorem by Kato on Naiver-Stokes equations}.
Revista Matem\"{a}tica Iber., 13(1997), 515-541.

\bibitem{Can2} M. Cannone,
{\it Harmonic analysis tools for solving the incompressible Naiver-Stokes equations},
Handbook of Mathematical fluid dynamics, Vol. 3,  Elsevier B.V., 2004, 161-244.

\bibitem{CMP}
M. Cannone, Y. Meyer and F. Planchon,
{\it Solutions autosimilaires des {\'e}quations de Navier-Stokes}.
{S{\'e}minaire ``{\'E}quations aux D{\'e}riv{\'e}es Partielles" de l'{\'E}cole polytechnique},
Expos{\'e} VIII, 1993--1994.

\bibitem{Chemin-book} J.-Y. Chemin,
{\it Perfect incompressible fluids}. Oxford University Press, New York, 1998.

\bibitem{CDGG1} J.-Y. Chemin, B. Desjardins, I. Gallagher and E. Grenier,
{\it Fluids with anisotropic viscosity. Special issue for R. Temam's 60th birthday}.
 M2AN Math. Model. Numer. Anal., 34(2000), no. 315--335.

\bibitem{CDGG2} J.-Y. Chemin, B. Desjardins, I. Gallagher and E. Grenier,
{\it Mathematical geophysics. An introduction to rotating fluids and the Navier-Stokes equations}.
Oxford University Press, New York, 2006.

\bibitem{CG}
J.-Y. Chemin and I. Gallagher, {\it On  the global wellposedness
of the 3-D Navier-Stokes equations with large initial data}.
Ann. de l'Ecole Norm. Sup., 39(2006), 679--698.

\bibitem{CZ} J.-Y. Chemin and P. Zhang, {\it On the global  wellposedness of the 3-D incompressible anisotropic
Navier-Stokes equations}. Comm. Math. Phys., 272(2007), 529--566.



\bibitem{CMZ} Q. Chen, C. Miao and Z. Zhang, {\it
Global well-posedness for the compressible Navier-Stokes equations
with the highly oscillating initial velocity}. To appear in Comm. Pure Appl. Math., 2009.

\bibitem{CMZ1} Q. Chen, C. Miao and Z. Zhang, {\it On the well-posedness for the viscous
 shallow water equations}, SIAM J. Math. Anal., 40(2008),443¨C--4.



\bibitem{CMZ2} Q. Chen, C. Miao and Z. Zhang, {\it
Well-posedness in critical spaces for the compressible Navier-Stokes equations with
density dependent viscosities}, To appear in Revista Matem\'atica Iberoamericana, 2009.


\bibitem{Fuj-Kat} H. Fujita and T. Kato, {\it On the Navier-Stokes initial value problem I}.
Arch. Rational Mech. Anal., 16(1964), 269-315.

\bibitem{Giga1} Y. Giga, K. Inui, A. Mahalov and S. Matsui, {\it Navier-Stokes equations in
a rotating frame in $\R^3$ with initial data nondecreasing at infity}. Hokkaido Math. J.,
35(2006), 321-364.

\bibitem{Giga2} Y. Giga, K. Inui, A. Mahalov, S. Matsui and J. Saal, {\it Rotating NS-equations
in $\R^3_+$ with initial data nondecreasing at infity: the Ekman boundary layer problem}.
Arch. Ration. Mech. Anal., 186(2007), 177--224.

\bibitem{Giga3} Y. Giga, K. Inui, A. Mahalov and J. Saal, {\it Uniform global solvability
of the rotating  Navier-Stokes equations for nondecaying initial data}.
Indiana Univ. Math. J., 57(2008), 2775--2791.


\bibitem{HS} M. Hieber and Y. Shibata, {\it The Fujita-Kato approach to the Navier-Stokes
equations in the rotational framework}. Math. Z, in press, 2009.

\bibitem{Kato}  T. Kato, {\it Strong $L^p$-solutions of Navier-Stokes equations in $\R^n$
with applications to weak solutions}. Math. Z, 187(1984), 471-480.

\bibitem{Majda} A. Majda, {\it An introducution to PDESs and Waves for the Atmosphere and Ocean}.
Courant Lecture Notes in Math., 2003.


\bibitem{OSA} D. Oliver, P. Stefanie  and V. Alexander, {\it A rotation method which gives linear $L^p$-estimates
for powers of the Ahlfors-Beurling operator}.
J. Math. Pures Appl., 86(2006), 492-509.



\bibitem{Paicu} M. Paicu, {\it \'{E}tude asymptotique pour les fluides anisotropes en rotation rapide dans le
cas p\'{e}riodique}.  J. Math. Pures Appl., 83(2004), 163--242.

\bibitem{Ped} J. Pedlosky, {\it Geophysical fluid dynamics}. Springer-Verlag, New York, 1987.




\end{thebibliography}
\end{document}